\documentclass[oneside]{amsart}
\usepackage[backref]{hyperref}
\usepackage{amsmath,ascmac,color,graphicx,amsfonts,amssymb,latexsym,amsthm,enumerate,comment,bm,pifont,mathtools}
\graphicspath{{./img_pdf/}}
\usepackage[hang,small,bf]{caption}
\usepackage[subrefformat=parens]{subcaption}
\usepackage[all]{xy} 
\usepackage{tikz}
\usetikzlibrary{intersections,calc,arrows.meta}
\usepackage{tikz}
\tikzset{
  point/.style={circle, fill=black, inner sep=1pt, minimum size=3pt}
}

\def\Z{\mathbb{Z} } 
\def\K{\mathbb{K} } 
 
\def\R{\mathbb{R} }

\def\snc{\mathrm{snc} }

\DeclareMathOperator{\rank}{rank}
\DeclareMathOperator{\im}{im}

\theoremstyle{definition}
\newtheorem{definition}{Definition}
\newtheorem{theorem}{Theorem}[section]
\newtheorem{proposition}[theorem]{Proposition}
\newtheorem{lemma}[theorem]{Lemma}
\newtheorem{corollary}[theorem]{Corollary}
\newtheorem{remark}{Remark}
\newtheorem{claim}{Claim}
\newtheorem{example}{Example}
\title{On Discrete Morse-Bott Theory}

\author{Yuto Nishikawa}

\author{Tomoo Yokoyama}

\date{\today}

\address{Department of Mathematics, Faculty of Science, Saitama University, Shimo-Okubo 255, Sakura-ku, Saitama-shi, 338-8570 Japan\\}

\date{\today}

\makeatletter
\@namedef{subjclassname@2020}{%
  \textup{2020} Mathematics Subject Classification}
\makeatother

\subjclass[2020]{Primary 37B30; Secondary 05C70}

\keywords{Morse-Bott function, discrete Morse theory, CW complexes, Poincar{\'e} polynomial}

\begin{document}
\maketitle

\begin{abstract}
This paper shows that discrete Morse-Bott theory can be developed as a natural extension of R. Forman's discrete Morse theory by improving the definition of the discrete Morse-Bott function introduced by S. Yaptieu. To this end, we demonstrate that the combinatorial structure of critical cells can be extended to critical sets intuitively. Furthermore, we establish the discrete Morse-Bott inequalities, providing a unified view that extends both the discrete Morse inequalities and the continuous Morse-Bott inequalities.
\end{abstract}

\section{Introduction}
Morse theory has been developed as a framework for describing the topological structure of manifolds via the critical points of smooth functions (cf.~\cite{milnor1963morse}). It establishes a correspondence between the properties of critical points of Morse functions and the topology of manifolds. R. Bott \cite{bott1954MBT} generalized Morse theory to Morse-Bott theory by allowing isolated critical points to be replaced by critical submanifolds. Morse-Bott theory is useful for studying functions with symmetries or conserved quantities and has applications in symplectic geometry (cf. \cite{austin1995morse,banyaga2004lectures}).

R. Forman \cite{forman1998DMT} introduced discrete Morse theory, which is a combinatorial version of Morse theory defined on CW complexes. Thus, it is suitable for computing homology and simplifying cell complexes, making it useful in various fields such as topological data analysis \cite{harker2014discrete,mischaikow2013morse}.

S. Yaptieu \cite{yaptieu2017discrete} proposed a generalized version of discrete Morse theory, called a discrete Morse-Bott function, where isolated critical points are replaced by more general subsets. This generalization reflects the structure of Morse-Bott theory in a discrete setting and allows the study of discrete functions with symmetry. 
However, although this work remains a preprint from 2017, the definition slightly differs from intuitive expectations, which implies the existence of some gaps (e.g. Remark~2.1(ii), the proofs of Lemmas~2.3, 2.4, 2.8). Therefore, we improve the definition to fill these gaps and make it more intuitive by extending the concept of critical points in discrete Morse theory to allow critical sets. 
In fact, we demonstrate a discrete Morse-Bott inequality, which naturally generalizes the discrete Morse inequalities known from the discrete Morse theory \cite{forman1998combinatorial,forman1998DMT} and can be derived from the continuous Morse-Bott theory \cite{David2008MBIneq}. The inequality gives an upper bound for the Betti numbers of the underlying complex based on the structure of the critical sets. Furthermore, it enables analysis of data using group actions or symmetries.

The present paper consists of four sections.
In the next section, we review the fundamental concepts and properties of discrete Morse theory. 
In \S~3, we introduce a discrete Morse-Bott function, which is a generalization of a discrete Morse function (Theorem~\ref{thmmorse-mb}). 
In the final section, we observe that various properties that hold for discrete Morse functions also hold for discrete Morse-Bott functions, and establish the discrete Morse-Bott inequality (Theorem~\ref{thmmbpoly}), which can be regarded as an extension of the discrete Morse inequalities (Proposition~\ref{rmkmb}).

\section{Preliminaries}

In this section, we review the fundamental concepts and properties of discrete Morse theory. 

\subsection{Fundamental concepts}

We define CW complex as follows. 

\begin{definition}
Let $X, Y$ be topological spaces and $f \colon A\to X$ a continuous map on a subspace $A \subseteq Y$. Then we define a subset 
\[
X\cup_f Y := (X\sqcup Y)/\sim,
\]
where $\sqcup$ denotes a disjoint union and the equivalence relation $\sim$ is the smallest one such that for every $a\in A$, we have $a\sim f(a)$. The map $f$ is called an \textbf{attaching map}. 
\end{definition}

\begin{definition}
A \textbf{CW complex} $X$ is constructed inductively as follows: \\
\textrm{(i)} $X^{-1}=\emptyset$. \\
\textrm{(ii)} For any integer $n \ge 0$, suppose the topological space $X^{n-1}$ has been defined. Then we set
\[
X^n := X^{n-1}\cup_{f^n}\coprod_{\lambda\in\Lambda^n}D^n_\lambda
\]
where $D^n$ is the $n$-dimensional closed disk and $f^n:\coprod_{\lambda\in\Lambda^n}\partial D_\lambda^n\to X^{n-1}$ is the attaching map. \\
We define $X := \bigcup_{n\in\Z_{\geq 0}} X^n$ and equip $X$ with the topology such that
\[
U\subset X \text{ is open if, for any } n\in\Z_{\geq 0}, \; U\cap X^n \text{ is open in } X^n.
\]
The map
\[
\varphi: D^n \hookrightarrow \coprod_{\lambda\in\Lambda^n}D^n_\lambda \twoheadrightarrow X^n \hookrightarrow X
\]
is called a \textbf{characteristic map} and $e^n := \varphi(\mathrm{Int}\,D^n)$ is called a \textbf{cell}, where $\hookrightarrow$ denotes an injection and $\twoheadrightarrow$ denotes a surjection. 
The integer $n$ is the dimension of the cell, denoted $\dim e^n$. \end{definition}

By definition, note that we have $X=\coprod_{\lambda\in\Lambda}e_\lambda^n$.
For simplicity, we identify a CW complex $\K$ with the family of all its cells.
With a slight abuse of notation, we also use $\K$ to denote this family.
From now on, let $\K$ be a CW complex.
We have the following concepts and notations. 

\begin{definition}
A cell $\sigma$ of $\K$ is a \textbf{face} of a cell $\tau$ of $\K$ (or $\sigma < \tau$) if $\sigma \subseteq \overline{\tau} - \tau$.  
\end{definition}

\begin{definition}
A face $\sigma$ of $\tau$ is a \textbf{facet} of $\tau$ (or $\sigma \prec \tau$) if $\,\dim \sigma=\dim \tau-1$.  
\end{definition}

\begin{definition}
A $k$-dimensional face $\sigma$ of $\tau$ is \textbf{regular} (or $\sigma \overset{\mathrm{reg}}{{<}} \tau$) if the following conditions hold: 
\\
\textrm{(1)} The restriction $\varphi_{\tau}\vert \colon\varphi^{-1}_{\tau}(\sigma)\,\rightarrow \, \sigma$ is a homeomorphism, where $D^{{\dim\tau}} \subset \R^{{\dim\tau}}$ is the closed unit disk and $\varphi_{\tau} \colon D^{{\dim\tau}} \to \K$ is the characteristic map of $\tau$. 
\\
\textrm{(2)} $\overline{\varphi^{-1}_{\tau}(\sigma) }$ is a $k$-dimensional closed disk.
\end{definition}

\begin{definition}
The relation $\sigma \overset{\mathrm{reg}}{{<}} \tau$ is denoted by $\sigma \overset{\mathrm{reg}}{{\prec}} \tau$ if $\dim \tau - \dim \sigma = 1$. 
\end{definition}

\begin{definition}
A cell $\sigma$ is an \textbf{irregular} face of $\tau$ (or $\sigma \overset{\mathrm{irr}}{{<}} \tau$) if $\sigma {<} \tau$ but $\sigma \overset{\mathrm{reg}}{\not<} \tau$.
\end{definition}

\begin{definition}
The relation $\sigma \overset{\mathrm{irr}}{{<}} \tau$ is denoted by $\sigma \overset{\mathrm{irr}}{{\prec}} \tau$ if $\dim \tau - \dim \sigma = 1$. 
\end{definition}

We provide the following examples to clarify the definitions.

\begin{example}
The following are examples of an irregular facet and a regular facet.
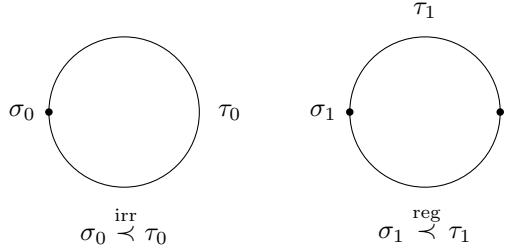
\[
\begin{tikzpicture}[point/.style={fill,
 shape=circle,inner sep=1pt,outer sep=0pt}]
\coordinate (A) at (1,0);
\draw (2,0) circle[radius=1];
\node[label=0 :$\tau_0$] at (3,0) {};
\node[point,label=180 :$\sigma_0$] at (A) {};
\node at (2,-1.5) {$\sigma_0\overset{\mathrm{irr}}{\prec}\tau_0$};

\coordinate (B) at (5,0);
\coordinate (C) at (7,0);
\draw (6,0) circle[radius=1];
\node[label=90 :$\tau_1$] at (6,1) {};
\node[point,label=180 :$\sigma_1$] at (B) {};
\node[point] at (C) {};
\node at (6,-1.5) {$\sigma_1\overset{\mathrm{reg}}{\prec}\tau_1$};

\end{tikzpicture}
\]
\end{example}

For a function on a CW complex, we call a cell noncritical when the value is not less than that on an adjacent higher-dimensional cell, as follows.
 
\begin{definition}
Let $f \colon \K \to \R$ be a function defined on a CW complex $\K$. 
A cell $\sigma$ is a \textbf{noncritical} facet of $\tau$ (or $\sigma \overset{\mathrm{nc}}{\prec} \tau$) if $\sigma \prec \tau\,$ and $\,f(\sigma)\geq f(\tau)$.
\end{definition}

We also introduce the following notations. 

\begin{definition}
For any cells $\nu<\tau$ (resp. $\nu \prec \tau$), when $f(\nu)>f(\tau)$, we write $\nu\overset{\mathrm{snc}}{<}\tau$ (resp. $\nu\overset{\mathrm{snc}}{\prec}\tau$). 
\end{definition}

Note that the concept of ``noncritical'' is also used in the definition of a discrete Morse-Bott function below.
We define the following ``interval''. 

\begin{definition}
For any cells $\sigma, \tau \in \K$ with $\sigma < \tau$, define the intervals $(\sigma, \tau) := \{ \alpha \in \K \mid \sigma < \alpha < \tau \}$, $[\sigma,\tau]:= \{ \alpha \in \K \mid \sigma \le \alpha \le \tau \}$, 
$[\sigma,\tau):= \{ \alpha \in \K \mid \sigma \le \alpha < \tau \}$, and 
$(\sigma,\tau]:= \{ \alpha \in \K \mid \sigma < \alpha \le \tau \}$.
\end{definition}

\subsection{Fundamental properties}

We have the following property. 

\begin{lemma}\label{lem:another_cell}
For any cells $\nu, \sigma, \tau \in \K$ with $\nu\overset{\mathrm{reg}}{\prec}\sigma\overset{\mathrm{reg}}{\prec}\tau$, there is a cell $\widetilde{\sigma}\ne\sigma$ such that $\nu\prec\widetilde{\sigma}\prec\tau$. 
\end{lemma}

The above lemma follows from the more general statement, Lemma~\ref{bnleqa} below, where a proof will be given.

\begin{lemma}\label{regface}
For any cells $\nu,\tau\in\K$ with $\nu \overset{\mathrm{reg}}{<}\tau$ and $\dim \tau - \dim \nu \geq 2$, there is a cell $\sigma \in (\nu,\tau)$ \textrm{(i.e.} $\sigma\in\K$ such that $\nu<\sigma<\tau$\textrm{)}. 
\end{lemma}

\begin{proof}
Assume that there are no such cells. 
Then $\varphi_{\tau}^{-1}(\nu)=\varphi_{\tau}^{-1}(\partial  \tau)$, which contradicts $\nu \overset{\mathrm{reg}}{<} \tau$. 
Therefore, there is a cell $\sigma \in \K$ such that $\nu<\sigma<\tau$. 

\end{proof}

\begin{lemma}\label{bnleqa}
For any cells $\nu,\alpha,\tau\in\K$ with $\nu<\alpha\overset{\mathrm{reg}}{\prec}\tau$, there is a cell $\beta \in (\nu,\tau)$ with $\beta\nleq\alpha$. 

\end{lemma}

\begin{proof}
By $\varphi_{\tau}^{-1}(\nu)\subset\varphi_{\tau}^{-1}(\overline{\alpha}) = \overline{\varphi_{\tau}^{-1}(\alpha)}$, we have $\varphi_{\tau}^{-1}(\nu) \subset \partial\varphi_{\tau}^{-1}(\alpha)$. 
From definition of ``regular'', there is a cell $\beta \in \K$ such that $\tau>\beta>\nu$ and $\beta\nleq\alpha$. 

\end{proof}

\begin{lemma}\label{nu<irrtau}
Let $\nu,\tau\in\K$ be cells with $\nu<\tau$. 
If there is exactly one cell $\alpha\in (\nu,\tau)$, 
then $\alpha\overset{\mathrm{irr}}{<}\tau$. 

\end{lemma}

\begin{proof}
Assume that $\alpha \overset{\mathrm{reg}}{<} \tau$. 
If $\dim \tau - \dim \alpha \geq 2$, then Lemma~\ref{regface} implies that there is a cell $\beta\in\K$ such that $\alpha<\beta<\tau$, which contradicts the uniqueness of $\alpha$. 
Thus $\dim \tau - \dim \alpha = 1$. 
However, Lemma~\ref{bnleqa} implies the existence of another cell between $\nu$ and $\tau$, which contradicts the uniqueness of $\alpha$. 
\end{proof}

\begin{remark}
In the previous lemma, the statement $\nu\overset{\mathrm{irr}}{\prec}\alpha$ need not hold. 
In fact, there is a CW complex that satisfies $X=\nu\cup\widetilde{\nu}\cup\alpha\cup\tau$, $\nu\overset{\mathrm{reg}}{\prec}\alpha$ as in Figure~\ref{fig:c_ex_01}. 

\begin{figure}[htbp]
\[
\begin{tikzpicture}
\coordinate (D) at (-210:2);
\coordinate (E) at (30:2);
\draw (D) --(E);
\node[point,label=180 : $\widetilde{\nu}$] at (D) {};
\node[point,label=0 :$\nu$] at (E) {};
\node[label=90 :$\alpha$] at (90:1) {};
\node[label=-90 :$\tau$] at (-90:2) {};
\draw (-210:2) arc[radius=2, start angle=-210,end angle=30];
\draw (180:2) arc[x radius=2,y radius=0.3, start angle=-180,end angle=0];
\draw[dash pattern=on 3pt off 2pt] (0:2) arc[x radius=2,y radius=0.3, start angle=0,end angle=180];
\end{tikzpicture}
\]
 \caption{An example of a CW complex without an irregular cell of one higher dimension for any zero cell}
  \label{fig:c_ex_01}
\end{figure}
\end{remark}

\subsection{Discrete Morse function}

In this subsection, we review the Discrete Morse function, which will be extended in the next section.

Let $f \colon \K \to \R$ be a function defined on a CW complex $\K$.
For any $k$-dimensional cell $\sigma^{(k)} \in \K$, set numbers $U(\sigma)$ and $D(\sigma)$ as follows: 
\[
\begin{split}
U(\sigma) := \# \{ \tau^{(k+1)} \mid \sigma \overset{\mathrm{nc}}{\prec} \tau \} = 
& \# \{ \tau^{(k+1)} \mid \sigma \prec \tau, f(\sigma) \geq f(\tau)\} 
\end{split}
\]
\[
D(\sigma):= \# \{ \nu^{(k-1)} \mid \nu \overset{\mathrm{nc}}{\prec} \sigma \} = \# \{ \nu^{(k-1)} \mid \nu \prec \sigma, f(\nu) \geq f(\sigma)\} 
\]

\begin{definition}[Discrete Morse function]\label{dmmorsefnc} 
The function $f$ is \textbf{discrete Morse} \cite{forman1998DMT} if the following conditions hold for any $\sigma \in \K$:
\\
\textrm{(M1)} For any cell $\tau$ with $\sigma \overset{\mathrm{irr}}{<} \tau$, we have 
$f(\sigma)<f(\tau)$. 
 \\
 \textrm{(M2)} $U(\sigma) \leq 1$.  
\\
 \textrm{(M3)} For any cell $\nu$ with $\nu \overset{\mathrm{irr}}{<} \sigma$, we have 
 $f(\nu)<f(\sigma)$. 
 \\
 \textrm{(M4)} $D(\sigma) \leq 1$.
\end{definition}

By definition, we have the following observations.

\begin{remark}
From definition (M1) of discrete Morse function, it follows that, for any discrete Morse function $f \colon \K \to \R$ and any cells $\sigma, \tau \in \K$, the condition $\sigma \overset{\mathrm{nc}}{\prec} \tau$ implies the condition $\sigma \overset{\mathrm{reg}}{\prec} \tau$.
\end{remark}

\begin{lemma}[Lemma~6.11.\cite{knudson2015morse}]\label{lemma:one}
For any discrete Morse function $f \colon \K \to \ \R$ and for any cell $\sigma \in \K$, we have the following inequality
\[
U(\sigma) + D(\sigma) \leq 1
\]
\end{lemma}

\subsection{Fundamental concepts}

We define fundamental concepts as follows. 

\begin{definition}
Let $f$ be a function defined on a CW complex. 
A cell $\sigma$ is \textbf{critical} if $U(\sigma) = D(\sigma) = 0$. 
\end{definition}

We recall the combinatorial vector field introduced by Forman \cite{forman1998combinatorial} as follows.

\begin{definition}
A map $V \colon \K \to \K \sqcup \{ 0\}$ on a CW complex $\K$ is a \textbf{combinatorial vector field} on $\K$ if it satisfies the following conditions: 
\\
\textrm{(1)} For any cell $\sigma \in V(\K) \setminus \{ 0\}$, we have $V(\sigma) = 0$. 
\\
\textrm{(2)} For any cell $\sigma \in \K - V^{-1}(0)$, we have $\sigma \overset{\mathrm{reg}}{\prec} V(\sigma)$.
\\
\textrm{(3)} For any cell $\sigma \in \K$, we have $\#  V^{-1}(\sigma) \leq 1$. 
\end{definition}

\begin{definition}
Let $V$ be a combinatorial vector field on a CW complex $\K$. 
Define arrows $\rightarrow$ as follows: 
$\sigma \rightarrow \tau$ if $V(\sigma)=\tau$. 
\end{definition}

Notice that $\sigma \rightarrow \tau$ implies $\sigma \overset{\mathrm{reg}}{\prec} \tau$. 

\begin{definition}
A path $ \sigma_0 \rightarrow \tau_0 \succ \sigma_1 \rightarrow \tau_1 \succ \dots \rightarrow \tau_{m-1} \succ \sigma_m$ is a \textbf{V-path} if for each $i$, $\sigma_i\ne\sigma_{i+1}$
\end{definition}

\begin{definition}
A V-path $ \sigma_0 \rightarrow \tau_0 \succ \sigma_1 \rightarrow \tau_1 \succ \dots \sigma_m \rightarrow \tau_m \succ \sigma_{m+1}$ is a \textbf{closed orbit} if $\sigma_{m+1} = \sigma_0$.
\end{definition}


\section{Discrete Morse-Bott function}

In this section, we introduce a new definition of discrete Morse-Bott function by improving the original definition so that it also allows comparisons between cells that differ by at least two dimensions. 
With this improvement, results that one would naturally expect as analogues of the continuous case can be shown in a precise way.

Therefore, we modify the definition of a discrete Morse-Bott function as follows. 

\subsection{Definition of discrete Morse-Bott function}
Let $f \colon \K \to \R$ be a function on a CW complex $\K$. 

For any $k$-dimensional cell $\sigma^{(k)}$, set $U^{\mathrm{snc}}(\sigma)$ and $D^{\mathrm{snc}}(\sigma)$ as follows: 
\[
\begin{split}
U^{\mathrm{snc}}(\sigma) := \# \{ \tau^{(k+1)} \mid \sigma \overset{\mathrm{snc}}{\prec} \tau \} = \# \{ \tau^{(k+1)} \mid \sigma \prec \tau, f(\sigma) > f(\tau)\} 
\end{split}
\]
\[
\begin{split}
D^{\mathrm{snc}}(\sigma):= \# \{ \nu^{(k-1)} \mid \nu \overset{\mathrm{snc}}{\prec} \sigma \} = \# \{ \nu^{(k-1)}  \mid \nu \prec \sigma, f(\nu) > f(\sigma)\} 
\end{split}
\]

Using these notations, we introduce a discrete Morse-Bott function as follows. 

\begin{definition}[Discrete Morse-Bott function]\label{dmbtfnc} 
A function $f \colon \K \to \R$ on a CW complex $\K$ is a \textbf{discrete Morse-Bott} function if the following conditions hold for any $\sigma^{(k)} \in \K$:
\\
\textrm{(M1)} For any cell $\tau$ with $\sigma \overset{\mathrm{irr}}{<} \tau$, we have 
$f(\sigma)<f(\tau)$. 
 \\
 \textrm{(MB2)} $U^{\mathrm{snc}}(\sigma) \leq 1$.  
\\
 \textrm{(M3)} For any cell $\nu$ with $\nu \overset{\mathrm{irr}}{<} \sigma$, we have 
 $f(\nu)<f(\sigma)$. 
 \\
 \textrm{(MB4)} $D^{\mathrm{snc}}(\sigma) \leq 1$.
\end{definition}

\subsubsection{Collections and paths}
 
\begin{definition}
For any $r \in \R$, a sequence $(\sigma_0, \sigma_1, \ldots , \sigma_l)$ is an \textbf{$\bm{r}$-path} from a cell $\sigma_0 \in \K$ to a cell $\sigma_l \in \K$ if $\{\sigma_0, \sigma_1, \ldots , \sigma_l\} \subseteq f^{-1}(r)$ and either $\sigma_i \prec \sigma_{i+1}$ or $\sigma_i \succ \sigma_{i+1}$ for any $i \in \{0,1, \ldots , l-1\}$. 
\end{definition}

Roughly speaking, an $r$-path is a path from a starting point to an endpoint that passes through cells of equal value whose dimensions differ by one, representing a kind of connectedness.

\begin{definition}[Collection]
 A \textbf{collection} $L$ for $f$ is a maximal family of cells satisfying the following conditions: 
\\
\textrm{(1)} There is a number $c \in \R$ such that $L \subset f^{-1}(c)$.
\\
\textrm{(2)} For any $\sigma, \tau \in L$, there is a $c$-path from $\sigma$ to $\tau$.
\end{definition}

Notice that any collection for $f$ coincides with an equivalence class for the equivalence relation given by the concept of $r$-path for some $r \in \R$.
Denote by $\mathcal{L}_f$ the set of collections of $f$. 
Then $\K = \bigsqcup \mathcal{L}_f$. 
For any cell $\sigma \in \K$, denote by $\mathcal{L}_f(\sigma)$ the collection containing $\sigma$. 

Notice that condition (2) in the previous definition is used only in Remark~\ref{rmk:3},  Theorem~\ref{thmmorse-mb}, and Proposition~\ref{rmkmb}.

For any collection $L$ for $f$ and any $k$-dimensional cell $\sigma^{(k)} \in L$, the following relations hold : 
\[
\begin{split}
U^{\mathrm{snc}}(\sigma) = \# \{ \tau^{(k+1)} \notin L \mid \sigma \overset{\mathrm{nc}}{\prec} \tau \} 
& = \# \{ \tau^{(k+1)} \notin L \mid \sigma \prec \tau, f(\sigma) \geq f(\tau)\}
\end{split}
\]
\[
\begin{split}
D^{\mathrm{snc}}(\sigma)= \# \{ \nu^{(k-1)} \notin L \mid \nu \overset{\mathrm{nc}}{\prec} \sigma \} & = \# \{ \nu^{(k-1)} \notin L \mid \nu \prec \sigma, f(\nu) \geq f(\sigma)\} 
\end{split}
\]

\tikzset{arrow/.style={->,thick}}
\tikzset{point/.style={fill,
shape=circle,inner sep=1pt,outer sep=0pt}}

Note that we often use symbols $\nu$, $\sigma$, and $\tau$ with the relation $\nu \prec \sigma \prec \tau$ among their dimensions, unless otherwise stated.

\begin{example}
The function on the left in Figure~\ref{fig:ex_DMF} is a discrete Morse-Bott function. 
On the other hand, the function on the right in Figure~\ref{fig:ex_DMF} is not, since it fails to satisfy (MB4) at a one-dimensional cell $\sigma$. 
Indeed, we obtain $U^{\mathrm{snc}}(\sigma)=2$, which contradicts (M4).

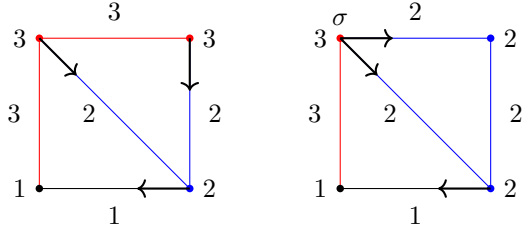
\begin{figure}
\[
\begin{tikzpicture}
\coordinate (A) at (0,0);
\coordinate (B) at (2,0);
\coordinate (C) at (0,2);
\coordinate (D) at (2,2);
\draw (A) --(B);
\draw[red] (A) --(C);
\draw[blue] (B) --(C);
\draw[blue] (B) --(D);
\draw[red] (C) --(D);
\node[point,label=180 :$1$] at (A) {};
\node[point,color=blue,label=0 :$2$] at (B) {};
\node[point,color=red,label=180 :$3$] at (C) {};
\node[point,color=red,label=0 :$3$] at (D) {};
\node[label=270 :$1$] at (1,0) {};
\node[label=180 :$3$] at (0,1) {};
\node[label=180 :$2$] at (1,1) {};
\node[label=0 :$2$] at (2,1) {};
\node[label=90 :$3$] at (1,2) {};
\draw[arrow] (0,2) --++(-45:0.7);
\draw[arrow] (2,2) --++(-90:0.7);
\draw[arrow] (2,0) --++(180:0.7);

\coordinate (A) at (4,0);
\coordinate (B) at (6,0);
\coordinate (C) at (4,2);
\coordinate (D) at (6,2);
\draw (A) --(B);
\draw[red] (A) --(C);
\draw[blue] (B) --(C);
\draw[blue] (B) --(D);
\draw[blue] (C) --(D);
\node[point,label=180 :$1$] at (A) {};
\node[point,color=blue,label=0 :$2$] at (B) {};
\node[point,color=red,label=180 :$3$,label=90 :$\sigma$] at (C) {};
\node[point,color=blue,label=0 :$2$] at (D) {};
\node[label=270 :$1$] at (5,0) {};
\node[label=180 :$3$] at (4,1) {};
\node[label=180 :$2$] at (5,1) {};
\node[label=0 :$2$] at (6,1) {};
\node[label=90 :$2$] at (5,2) {};
\draw[arrow] (4,2) --++(0:0.7);
\draw[arrow] (4,2) --++(-45:0.7);
\draw[arrow] (6,0) --++(180:0.7);
\end{tikzpicture}
\]
\caption{Left: discrete Morse-Bott function; right: non-discrete Morse-Bott function.}
\label{fig:ex_DMF}
\end{figure}

\end{example}

We have the following observation.

\begin{lemma}\label{dmfisdmbf}
Let $f \colon \K \to \R$ be a discrete Morse function. 
The function $f$ is a discrete Morse-Bott function such that $\#  L\le 2$ for any collection $L\in\mathcal{L}_f$. 
\end{lemma}

\begin{proof}
Fix any $\sigma\in\K$. 
conditions (M1) and (M3) follow from definition of a discrete Morse function. 
The following conditions 
\[
U^\snc(\sigma)\le U(\sigma)\le1
\]
\[
D^\snc(\sigma)\le D(\sigma)\le1
\]
imply conditions (MB2) and (MB4). 
This means that $f$ is a discrete Morse-Bott function.

Assume that $\#  L\ge 3$. 
Then there is a $c$-path $\sigma_0,\sigma_1,\sigma_2$. 
Therefore, we obtain that $U(\sigma_1)+D(\sigma_1)\ge 2$, which contradicts $U(\sigma_1)+D(\sigma_1) = 1$ because of Lemma~\ref{lemma:one}. 
Thus $\#  L\le 2$. 
\end{proof}

\subsection{Reduced collection}

Let $f \colon \K \to \R$ be a discrete Morse-Bott function. 
Recall that a cell $\sigma \in \K$ is \textbf{critical} if and only if $U(\sigma)=D(\sigma)=0$. 
We introduce the following concepts. 

\begin{definition}
A cell $\sigma \in \K$ is \textbf{upward noncritical} if $U^\snc(\sigma)=1$. 
\end{definition}

\begin{definition}
A cell $\sigma \in \K$ is \textbf{downward noncritical} if $D^\snc(\sigma)=1$. 
\end{definition}


\begin{definition}
A cell $\sigma \in \K$ is \textbf{weakly critical} if $D^\snc(\sigma)=U^\snc(\sigma)=0$.
\end{definition}

\begin{definition}[Reduced collection]

 A \textbf{reduced collection} $C$ for $f$ is a maximal family of cells satisfying the following conditions: 
\\
\textrm{(1)} There is a number $c \in \R$ such that $C \subset f^{-1}(c)$.
\\
\textrm{(2)} For any $\sigma, \tau \in C$, there is a $c$-path from $\sigma$ to $\tau$.
\\
\textrm{(3)} For any $\sigma\in C$, $\sigma$ is weakly critical.
\end{definition}

Denote by $\mathcal{C}_f$ the set of reduced collections of $f$.
For any cell $\sigma \in \K$, denote by $\mathcal{C}_f(\sigma)$ the reduced collection containing $\sigma$.

\begin{remark}
For any reduced collection $C\in\mathcal{C}_f$, there exists a unique collection $L$ such that $C\subset L$.
\end{remark}

Denote by $L(C)$ the collection containing $C$.

\begin{example}
We describe below two examples of reduced collections.
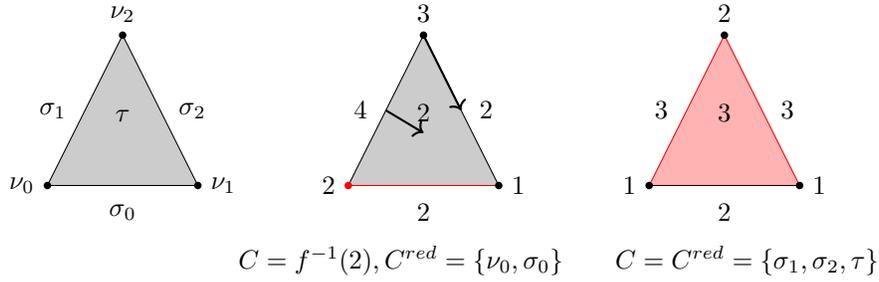
\[
\begin{tikzpicture}
\coordinate (G) at (-4,0);
\coordinate (H) at (-2,0);
\coordinate (I) at (-3,2);
\draw (G) --(H) --(I) --cycle;
\fill[opacity=0.2] (G) --(H) --(I) --cycle;
\node[point,label=180 :$\nu_0$] at (G) {};
\node[point,label=0 :$\nu_1$] at (H) {};
\node[point,label=90 :$\nu_2$] at (I) {};
\node[label=180 :$\sigma_1$] at (-3.5,1) {};
\node[label=0 :$\sigma_2$] at (-2.5,1) {};
\node[label=-90 :$\sigma_0$] at (-3,0) {};
\node[label=90 :$\tau$] at (-3,0.6) {};

\coordinate (A) at (0,0);
\coordinate (B) at (2,0);
\coordinate (C) at (1,2);
\draw[red] (A) --(B);
\draw (B) --(C);
\draw (A) --(C);
\fill[opacity=0.2] (A) --(B) --(C) --cycle;
\node[point,color=red,label=180 :$2$] at (A) {};
\node[point,label=0 :$1$] at (B) {};
\node[point,label=90 :$3$] at (C) {};
\node[label=180 :$4$] at (0.5,1) {};
\node[label=0 :$2$] at (1.5,1) {};
\node[label=-90 :$2$] at (1,0) {};
\node[label=90 :$2$] at (1,0.6) {};
\draw[arrow] (1,2) --(1.5,1);
\draw[arrow] (0.5,1) --(1,0.7);
\node at (0.7,-1) {$L=f^{-1}(2),C=\{\nu_0, \sigma_0\}$};

\coordinate (D) at (4,0);
\coordinate (E) at (6,0);
\coordinate (F) at (5,2);
\draw (D) --(E);
\draw[red] (F) --(E);
\draw[red] (D) --(F);
\fill[opacity=0.3][red] (D) --(E) --(F) --cycle;
\node[point,label=180 :$1$] at (D) {};
\node[point,label=0 :$1$] at (E) {};
\node[point,label=90 :$2$] at (F) {};
\node[label=180 :$3$] at (4.5,1) {};
\node[label=0 :$3$] at (5.5,1) {};
\node[label=-90 :$2$] at (5,0) {};
\node[label=90 :$3$] at (5,0.6) {};
\node at (5.3,-1) {$L=C=\{\sigma_1, \sigma_2, \tau\}$};

\end{tikzpicture}
\]
\end{example}

\begin{remark}\label{rmk:3}
The following are equivalent for any discrete Morse-Bott function $f \colon \K \to \R$ and any cell $\sigma \in \K$:
\\
\textrm{(1)} $\sigma$ is critical.
\\
\textrm{(2)} $L(\mathcal{C}_f(\sigma))=\mathcal{C}_f(\sigma)$. 

\end{remark}

\subsection{Generalization of a discrete Morse function}

We have the following equivalence, which demonstrates that the discrete Morse-Bott function is a generalization of the discrete Morse function. 

\begin{theorem}\label{thmmorse-mb}
The following are equivalent for any  function $f \colon \K \to \R$: 
\\
\textrm{(1)} The function $f$ is discrete Morse. 
\\
\textrm{(2)} The function $f$ is discrete Morse-Bott such that, for any collection $L \in \mathcal{L}_f$, one of the following holds:
\\
\textrm{(i) } $\#L=1$.
\\
\textrm{(ii)} $\#L=2$ and $L$ is a reduced collection.
\end{theorem}

\begin{proof}
Suppose that $f$ is discrete Morse. 
By Proposition~\ref{dmfisdmbf}, we have that $\#  L\le 2$ for any collection $L\in\mathcal{L}_f$. 
If $\#  L = 2$, then Lemma~\ref{lemma:one} implies that $0= U^\snc(\sigma)+D^\snc(\sigma)$ and so that $L$ is a reduced collection, which implies assertion (2).

Conversely, suppose that assertion (2) holds. 
Fix any collection $L \in \mathcal{L}_f$ and any cell $\sigma \in L$. 
If $L$ is a singleton (i.e. $L = \{ \sigma \}$), then $U(\sigma)= U^\snc(\sigma)$ and $D(\sigma)= D^\snc(\sigma)$ and so (MB2) (resp. (MB4)) for $L$ implies (M2) (resp. (M4)) for $L$. 
If $\#  L =2$, then 
\begin{align*}
U(\sigma) =\# \{ \tau \mid \sigma \overset{\mathrm{nc}}{\prec} \tau \} &= \# \{ \tau\notin L \mid \sigma \overset{\mathrm{nc}}{\prec} \tau \}+\# \{ \tau\in L \mid \sigma \overset{\mathrm{nc}}{\prec} \tau \}\\
&\le0+1=1\\
D(\sigma) =\# \{ \nu \mid \nu \overset{\mathrm{nc}}{\prec} \sigma \} &= \# \{ \nu\notin L \mid \nu \overset{\mathrm{nc}}{\prec} \sigma \}+\# \{ \nu\in L \mid \nu \overset{\mathrm{nc}}{\prec} \sigma \}\\
&\le0+1=1
\end{align*}
and so (M2) and (M4) hold, because $\sigma$ is weakly critical. 
This means that $f$ is discrete Morse. 
\end{proof}

\begin{corollary}\label{cormorse-mb}
Every Morse function $f$ is a discrete Morse-Bott function such that, for any reduced collection $C\in\mathcal{C}_f$, $\#C\le2$ and $C=L(C)$. 
\end{corollary}

\begin{proof}
    By Theorem~\ref{thmmorse-mb}, we have that $\#C\le\#L(C)\le2$ for any reduced collection $C\in\mathcal{C}_f$. If $\#C=2$, then $\#L(C)=2$ and so $C=L(C)$. 
Thus, we may assume that $\#C=1$. 
If $\#L(C)>1$, then Theorem~\ref{thmmorse-mb} implies that $\#L(C)=2$ and $C=L(C)$, which contradicts the assumption $\#C=1$.
Therefore, $\#L(C)=1$ and so $C=L(C)$.
\end{proof}

\section{Natural properties of discrete Morse-Bott functions}

In this section, we observe that various properties that hold for discrete Morse functions also hold for discrete Morse-Bott functions.
To state them, we present the following technical lemmas.

\begin{lemma}\label{lemupnon}
Let $f \colon \K \to \R$ be a discrete Morse-Bott function and $\sigma\in \K$ an upward noncritical cell. Then any cell $\nu\in \K$ with $\nu\prec\sigma$ and $f(\nu)=f(\sigma)$ is upward noncritical. 

\end{lemma}

\begin{proof}
Since $\sigma\in \K$ is upward noncritical, there is a cell $\tau\in\K$ with $\tau\overset{\mathrm{snc}}{\succ}\sigma$. Therefore, we have $\nu\overset{\mathrm{reg}}{\prec}\sigma\overset{\mathrm{reg}}{\prec}\tau$. By Lemma~\ref{lem:another_cell}, there is a cell $\widetilde{\sigma}\ne\sigma$ such that $\nu\prec\widetilde{\sigma}\prec \tau$. Since $D^\snc(\tau)=1$, we have $\widetilde{\sigma}\overset{\mathrm{snc}}{\not\prec}\tau$. The following inequality holds: 
\[
f(\nu) = f(\sigma)> f(\tau) \ge f(\widetilde{\sigma})
\]
From $\nu\prec\widetilde{\sigma}$, the cell $\nu$ is upward noncritical.
\end{proof}

We have the following trichotomy. 

\begin{lemma}\label{exclude}
Let $f \colon \K \to \R$ be a discrete Morse-Bott function. 
Then one of the following statements holds exclusively for any cell $\sigma\in\K$: 
\\
\textrm{(1)} $\sigma$ is upward noncritical. 
\\
\textrm{(2)} $\sigma$ is downward noncritical.  
\\
\textrm{(3)} $\sigma$ is weakly critical.  
\end{lemma}

\begin{proof}
Suppose that $\sigma$ is not weakly critical. 
Assume that the cell $\sigma$ satisfies both (1) and (2). 
Then there are cells $\nu,\tau\in\K$ such that $\nu\overset{\mathrm{snc}}{\prec}\sigma\overset{\mathrm{snc}}{\prec}\tau$. 
Therefore, we have $\nu\overset{\mathrm{reg}}{\prec}\sigma\overset{\mathrm{reg}}{\prec}\tau$. 
By Lemma~\ref{lem:another_cell}, there is a cell $\widetilde{\sigma}\ne\sigma $ such that $\nu\prec\widetilde{\sigma}\prec\tau$. 
This implies that 
\[
f(\tau)<f(\sigma)<f(\nu)\le f(\widetilde{\sigma})\le f(\tau),
\]
which is a contradiction.  
\end{proof}

\begin{lemma}\label{nu>tau}
Let $f \colon \K \to \R$ be a discrete Morse-Bott function. 
For any cells $\nu,\tau\in\K$ with $\nu\overset{\mathrm{snc}}{<}\tau$, 
there is $\sigma\in (\nu,\tau]$ with $\nu \overset{\mathrm{snc}}{\prec} \sigma$.
\end{lemma}

\begin{proof}
Fix any cells $\nu,\tau\in\K$ with $\nu\overset{\mathrm{snc}}{<}\tau$. 
Then $f(\tau) < f(\nu)$.
Put $n := \dim \tau - \dim \nu$. 
We demonstrate the assertion for $\nu \overset{\mathrm{snc}}{<} \tau$ by induction on $n= \dim \tau - \dim \nu$.
If $ n= 1$, then one can choose $\sigma := \tau$. 
Thus, we may assume that $\dim \tau - \dim \nu > 1$ and the assertion for any pair $\nu' \overset{\mathrm{snc}}{<} \tau'$ with $\dim \tau' - \dim \nu' < n$ holds. 
By Lemma~\ref{regface}, there is a cell $\alpha\in (\nu,\tau)$. 
Write $m := \dim \tau - \dim \alpha < n$. 

Suppose that $m= 1$.
If $f(\alpha)\le f(\tau)$, then $f(\alpha)\le f(\tau) < f(\nu)$ and so 
the inductive hypothesis for $\nu \overset{\mathrm{snc}}{<} \alpha$ implies the assertion. 
Thus, we may assume that $f(\alpha) > f(\tau)$. 
Then $\alpha \overset{\mathrm{reg}}{\prec} \tau$.
By Lemma~\ref{bnleqa}, there is a cell $\beta \in (\nu,\tau)$ with $\beta\nleq\alpha$.
If $f(\beta) \le f(\tau)$, then $f(\beta) \le f(\tau) < f(\nu)$ and so the assertion holds by induction on $n$, because $\beta \in (\nu,\tau)$. 
Thus, we may assume that $f(\beta) > f(\tau)$. 
Then $\beta \overset{\mathrm{snc}}{<} \tau$.
By the inductive hypothesis for $\beta \overset{\mathrm{snc}}{<} \tau$, 
there is a cell $\beta' \in (\beta,\tau]\subset(\nu,\tau]$ such that $\beta \overset{\mathrm{snc}}{\prec} \beta'$. 
Suppose that $f(\beta') \leq f(\tau)$. 
Then $f(\beta') \leq f(\tau) < f(\nu)$ and so $\nu \overset{\mathrm{snc}}{<} \beta'$. 
The inductive hypothesis for $\nu \overset{\mathrm{snc}}{<} \beta'$ implies that there is $\sigma\in (\nu,\beta'] \subset (\nu,\tau]$ with $\nu \overset{\mathrm{snc}}{\prec} \sigma$.
Thus, we may assume that $f(\beta') > f(\tau)$. 
Then the inductive hypothesis for $\beta' \overset{\mathrm{snc}}{<} \tau$ implies that there is $\beta''\in (\beta',\tau]$ with $\beta' \overset{\mathrm{snc}}{\prec} \beta''$. Therefore, we have $\beta\overset{\mathrm{snc}}{\prec}\beta'\overset{\mathrm{snc}}{\prec}\beta''$, which contradicts Lemma~\ref{exclude}.

Suppose that $m>1$.
If $f(\alpha)\le f(\tau)$, then $f(\alpha) \leq f(\tau) < f(\nu)$
and $\nu < \alpha$, which implies the assertion from the inductive hypothesis for $\nu \overset{\mathrm{snc}}{<} \alpha$. 
Thus, we may assume that $f(\alpha)>f(\tau)$. 
Lemma~\ref{regface} implies that there is a cell $\widetilde{\alpha} \in (\alpha, \tau) \subset (\nu,\tau)$. 
The inductive hypothesis on $m$ implies the assertion holds, because $\dim\tau-\dim\widetilde{\alpha}<m$. 
\end{proof}

\begin{lemma}\label{unclem}
Let $f \colon \K \to \R$ be a discrete Morse-Bott function and $C \in \mathcal{C}_f$ a reduced collection. 
Then every cell $\sigma\in\overline{C}\setminus C$ with $f(\sigma) = f(C)$ is upward noncritical. 
\end{lemma}

\begin{proof}
By $\sigma\in \overline{C}$, there is a cell $\tau \in C$ such that $\sigma<\tau$. 
Choose such a cell $\tau_0\in C$ whose dimension is minimal. 
Since $f(\sigma) = f(C)$, we obtain $f(\sigma) = f(\tau_0)$. 

\begin{claim}
If $\dim \tau_0 - \dim \sigma = 1$, then $\sigma$ is upward noncritical.
\end{claim}
\begin{proof}
Suppose that $\dim \tau_0 - \dim \sigma = 1$.
Assume that $\sigma$ is downward noncritical. 
Then there is a cell $\nu\in\K$ with $\nu\overset{\mathrm{snc}}{\prec}\sigma$. 
By Lemma~\ref{bnleqa}, there is a cell $\widetilde{\sigma} \in (\nu,\tau_0)$ with $\widetilde{\sigma}\ne\sigma$. 
From $\tau_0\in C$, since $C$ is a reduced collection, we have $f(\widetilde{\sigma})\le f(\tau_0)$. 
Therefore, we have $f(\nu) \leq f(\widetilde{\sigma}) \leq f(\tau_0) = f(\sigma) < f(\nu)$, which is a contradiction. 
\end{proof}

By the previous claim, we may assume that $\dim \tau_0 - \dim \sigma > 1$. 

\begin{claim}
If there is no cell $\alpha \in  (\sigma,\tau_0)$ such that $\alpha \overset{\mathrm{reg}}{\prec} \tau_0$, then $\sigma$ is upward noncritical.
\end{claim}

\begin{proof}
Suppose that there is no cell $\alpha \in  (\sigma,\tau_0)$ such that $\alpha \overset{\mathrm{reg}}{\prec} \tau_0$.
By $\sigma \overset{\mathrm{reg}}{<} \tau_0$, there is a cell $\beta\in (\sigma,\tau_0)$ such that $\beta\overset{\mathrm{irr}}{<}\tau_0$. 
Then $f(\sigma) = f(\tau_0) > f(\beta)$. 
By $\sigma\overset{\mathrm{snc}}{<}\beta$, Lemma~\ref{nu>tau} implies that there is a cell $\widetilde{\sigma} \in (\sigma,\beta)$ such that $\sigma\overset{\mathrm{snc}}{\prec}\widetilde{\sigma}$. 
Therefore, $\sigma$ is upward noncritical.
\end{proof}

By the previous claim, we may assume that there is a cell $\alpha \in (\sigma,\tau_0)$ such that $\alpha \overset{\mathrm{reg}}{\prec}\tau_0$. 
From $\tau_0\in C$, we obtain $f(\alpha)\le f(\tau_0)$. 
Since $f(\tau_0) = f(\sigma)$, if $f(\alpha)<f(\tau_0)$, then $\sigma\overset{\mathrm{snc}}{<}\alpha$ and so Lemma \ref{nu>tau} implies that $\sigma$ is upward noncritical.
Thus, we may assume that $f(\alpha) = f(\tau_0)=f(C)$. 

\begin{claim}
The cell $\alpha$ is upward noncritical.
\end{claim}

\begin{proof}
Assume that $\alpha$ is downward noncritical.
Then there is a cell $\beta \overset{\mathrm{snc}}{\prec} \alpha$. 
By $\beta \prec \alpha \overset{\mathrm{reg}}{\prec} \tau_0$, Lemma~\ref{bnleqa} implies that there is a cell $\widetilde{\alpha} \neq \alpha$ with $\beta \prec \widetilde{\alpha} \prec \tau_0$. 
From $\tau_0 \in C$, we have the following contradiction: 
\[
f(\beta) > f(\alpha ) = f(\tau_0) \geq f(\widetilde{\alpha}) \geq f(\beta)
\]
\end{proof}

By the previous claim, there is a cell $\gamma$ with  $\alpha \overset{\mathrm{snc}}{\prec} \gamma$. 
Since $f(\sigma) = f(C) = f(\alpha) > f(\gamma)$, Lemma~\ref{nu>tau} implies that $\sigma$ is upward noncritical. 
\end{proof}

\begin{lemma}
Let $f$ be a discrete Morse-Bott function. 
For any reduced collection $C\in\mathcal{C}_f$, the set difference $\overline{C} \setminus C$ is a subcomplex of $\K$.
\end{lemma}

\begin{proof}
Put $A:=\overline{C}\setminus C$. 
Fix any cells $\nu<\sigma$ with $\sigma\in A$.
It suffices to show that $\nu\in A$. 
By definition of closure, we have $\nu\in\overline{C}$.
If $f(\nu)\ne f(C)$, then $\nu\notin C$ and so $\nu \in \overline{C}\setminus C = A$.
Thus, we may assume that $f(\nu)=f(C)$. 
By $\sigma\in \overline{C}\setminus C$, there is a cell $\tau \in C$ with $\sigma<\tau$. 
Choose such a cell $\tau_0 \in C$ whose dimension is minimal. 
If $f(\sigma)<f(\tau_0)=f(\nu)$, then the fact $\nu\overset{\snc}{<}\sigma$ and Lemma~\ref{nu>tau} imply $\nu\notin C$ and so $\nu\in A$. 
Suppose that $f(\sigma)=f(C)$.
By Lemma~\ref{unclem}, the cell $\sigma \in \overline{C}\setminus C$ is upward noncritical. 
There is a cell $\sigma' \in \K$ with $\sigma \prec \sigma'$ and $f(\sigma) > f(\sigma')$. 
Since $f(\sigma') < f(\sigma)=f(C) = f(\nu)$ and $\nu < \sigma \prec \sigma'$, we have $\nu \overset{\snc}{<} \sigma'$. 
Lemma~\ref{nu>tau} imply $\nu\notin C$ and so $\nu\in A$. 
Thus, we may assume that $f(\sigma)>f(C)$. 
Then $f(\sigma)>f(C) = f(\tau_0)=f(\nu)$. 

\begin{claim}
There is a cell $\alpha\in(\sigma,\tau_0)$ such that $f(\alpha)\le f(\tau_0)$ and $\alpha\notin C$.
\end{claim}

\begin{proof}
Assume that there is no cell $\alpha\in(\sigma,\tau_0)$ such that $f(\alpha)\le f(\tau_0)$ and $\alpha\notin C$.
By $\tau_0 \in C$, we have that $\sigma \overset{\mathrm{reg}}{<} \tau_0$ and $\dim \tau_0 - \dim \sigma \geq 2$. 
Applying Lemma~\ref{regface} to $\sigma \overset{\mathrm{reg}}{<} \tau_0$, there is a cell $\alpha_0 \in(\sigma,\tau_0)$. 
By $\tau_0 \in C$ and $\alpha_0 < \tau_0$, we obtain that $\alpha_0 \in A$, because $\tau_0>\sigma$ is a minimal cell which is in $C$.
By the assumption, we obtain that $f(\alpha_0) > f(\tau_0)$. 
Then $\alpha_0 \overset{\snc}{<} \tau_0$. 
Lemma~\ref{nu>tau} imply that there is a cell $\alpha_1 \in (\alpha_0, \tau_0) \subset (\sigma,\tau_0)$ with $\alpha_0 \overset{\snc}{<} \alpha_1$. 
By the assumption and the minimality of $\tau_0$, we obtain that $\alpha_1 \overset{\snc}{<} \tau_0$. 
By finite iterations of this process, one can obtain a cell $\alpha^+ \in (\sigma,\tau_0)$ with $\alpha^+ \overset{\snc}{\prec} \tau_0$, which contradicts that $\tau_0$ is weakly critical. 
\end{proof}

If $f(\alpha)< f(\tau_0)$, then the fact $\nu\overset{\snc}{<}\alpha$ and Lemma~\ref{nu>tau} imply $\nu\notin C$ and so $\nu\in A$. 
Thus, we may assume that $f(\alpha) = f(\tau_0)$. 
Then $f(\alpha) = f(\tau_0) = f(C) = f(\nu)$. 
By $\tau_0 \in C$ and $\alpha < \tau_0$, we obtain that $\alpha \in A$, because $\tau_0>\sigma$ is a minimal cell which is in $C$.
Since $f(\alpha)= f(C)$ and $\alpha\in A = \overline{C}\setminus C$, Lemma~\ref{unclem} implies that $\alpha \in \overline{C}\setminus C$ is upward noncritical. 
There is a cell $\alpha' \in \K$ with $\alpha \prec \alpha'$ and $f(\alpha) > f(\alpha')$. 
Since $f(\alpha') < f(\alpha)=f(C) = f(\nu)$ and $\nu < \alpha \prec \alpha'$, we have $\nu \overset{\snc}{<} \alpha'$. 
Lemma~\ref{nu>tau} imply $\nu\notin C$ and so $\nu\in A$. 
\end{proof}

For any cells $\tau, \sigma \in \K$, denote by  $[\tau:\sigma]$ the incidence number.  
For any reduced collection $C$ and any integer $k \in \Z_{\geq 0}$, denote by $C_k(C;\Z)$ the free $\Z$-module generated by the $k$-cells of $C$.  
Put $C_{-1}(C;\Z) := \{ 0 \}$.

\begin{definition}
For any reduced collection $C\in\mathcal{C}_f$ of a discrete Morse-Bott function $f$ and any integer $k \in \Z_{\geq 0}$, we define the boundary homomorphism  
\[
\partial_k^C:C_k(C;\Z)\rightarrow C_{k-1}(C;\Z)
\]
for any $k$-cell $\tau \in C$ by  
\[
\partial_k^C(\tau):=\sum_{\substack{\sigma\in C\\\sigma\prec\tau}}[\tau:\sigma]\sigma.
\]
\end{definition}

In the previous definition, the homomorphism is well-defined and independent of the orientation of the cells.  

\begin{proposition}
For any reduced collection $C\in\mathcal{C}_f$ of a discrete Morse-Bott function $f$ and any integer $k \in \Z_{>0}$, the equality 
\[
\partial_{k-1}^C \circ \partial_k^C = 0
\]
holds. 
\end{proposition}

\begin{proof}
Put
\[
X := \overline{C}, \quad 
A := \overline{C} \setminus C.
\] 
Then both $X$ and $A$ are subcomplexes of $\K$. 
Since the restriction of the boundary homomorphism 
\[
\partial_k : C_k(\K) \to C_{k-1}(\K)
\] 
to any subcomplex is also a boundary homomorphism, the following commutative diagram holds:
\[
  \xymatrix{
        & 0 \ar[d] & 0 \ar[d] & 0 \ar[d] &  \\
    \ar[r] & C_{k+1}(A) \ar[d] \ar[r] & C_{k}(A) \ar[r] \ar[d] & C_{k-1}(A) \ar[r] \ar[d] &  \\
    \ar[r] & C_{k+1}(X) \ar[d] \ar[r] & C_{k}(X) \ar[r] \ar[d] & C_{k-1}(X) \ar[r] \ar[d] &  \\
    \ar[r] & C_{k+1}(X,A) \ar[d] \ar[r] & C_{k}(X,A) \ar[r] \ar[d] & C_{k-1}(X,A) \ar[r] \ar[d] &  \\
        & 0 & 0 & 0 &   \\
  }
\]
Here, the columns in the above diagram are short exact sequences, and the rows are chain complexes.  
Moreover, the boundary homomorphism $\partial^C$ satisfies the following commutative diagram:
\[
  \xymatrix{
    \ar[r] & C_{k+1}(X,A) \ar[d] \ar[r] & C_{k}(X,A) \ar[r] \ar[d] & C_{k-1}(X,A) \ar[r] \ar[d] &  \\
    \ar[r] & C_{k+1}(C) \ar[r]^{\partial^C} & C_{k}(C) \ar[r]^{\partial^C} & C_{k-1}(C) \ar[r] &  \\
  }
\]
Since $X \setminus A = C$, we have an isomorphism 
\[
C_k(C) \cong C_k(X,A).
\] 
Therefore, it follows that 
\[
\partial^C \circ \partial^C = 0.
\]
\end{proof}

\begin{definition}
For any reduced collection $C\in\mathcal{C}_f$ of a discrete Morse-Bott function $f$, define the $k$-th Betti number $b_k^C$ of $C$ and the Poincar{\'e} polynomial $P_t(C)$ of $C$ as follows: 
\[
b_k^C := \rank\left(Z_k^C / B_k^C\right),
\]
\[
P_t(C) := \sum_k b_k^C t^k.
\] 
where $Z_k^C := \ker \partial_k^C$ and $B_k^C := \im \partial_{k+1}^C$. 
\end{definition}

\begin{definition}
For any discrete Morse function $f \colon \K \to \R$, define a combinatorial vector field $-\nabla f \colon \K \to \K \sqcup \{\emptyset\}$ as follows:  
\[
-\nabla f(\sigma) :=
\begin{cases}
\tau & \text{if } \sigma \overset{\mathrm{nc}}{\prec} \tau \\
\emptyset & \text{otherwise}
\end{cases}
\]
\end{definition}

\begin{definition}
For any discrete Morse-Bott function $f \colon \K \to \R$, define a combinatorial vector field $-\nabla_s f \colon \K \to \K \sqcup \{\emptyset\}$ as follows:  
\[
-\nabla_s f(\sigma) :=
\begin{cases}
\tau & \text{if } \sigma \overset{\mathrm{snc}}{\prec} \tau \\
\emptyset & \text{otherwise}
\end{cases}
\]
\end{definition}

\begin{definition}
Let $V$ be a combinatorial vector field.  
We say that $V$ is \textbf{positively bounded} if the following holds for any $\sigma \in \K$:
\[
\sup \left( \{ 0 \} \cup \left\{ r  \in \Z \, \middle| \,
\begin{aligned}
   &\text{there is a $V$-path containing no closed orbit}  \\
   &\text{of the form } \sigma = \sigma_0 \rightarrow \tau_0 \succ \sigma_1 \rightarrow 
     \dots \succ \sigma_r
\end{aligned}
\right\} \right) < \infty.
\]
\end{definition}

We have the following characterization, which is essentially demonstrated in \cite[p.95, Theorem~6.19]{knudson2015morse}. 

\begin{lemma}\label{lemgvf}
The following are equivalent for any positively bounded combinatorial vector field $V$:  
\\
\textrm{(1)} There is a discrete Morse function $f$ such that $V = -\nabla f$.  
\\
\textrm{(2)} $V$ has no closed orbits.  
\end{lemma}

We omit the proof of the previous lemma, because it proceeds in exactly the same way as the original proof of \cite[p.95, Theorem~6.19]{knudson2015morse}, except that the formula on line 19 
\[
f_p(\sigma^{p-1})=f_{p-1}(\sigma)+\frac{d(\sigma)}{2D+1}
\]
is replaced as follows: 
\[
f_p(\sigma^{p-1})=f_{p-1}(\sigma)+\sum_{k=1}^{d(\sigma)}\frac{1}{2^k}
\]

\begin{lemma}\label{lemsgvf}
Let $f$ be a discrete Morse-Bott function.  
If $-\nabla_s f$ is positively bounded, then there is a discrete Morse function $g$ with $
-\nabla_s f = -\nabla g$.
\end{lemma}

\begin{proof}
Since $-\nabla_s f$ has no closed orbits, Lemma~\ref{lemgvf} implies the assertion.
\end{proof}

We introduce the following notation to clarify the following discussion.

\begin{definition}
Let $f:\K\to\R$ be a function, $k \in \{0,1,\dots,\dim \K\}$ an integer, and $J$ a subcomplex of $\K$. 
Then $c^J_k,m^J_k,u^J_k,d^J_k$ are defined as follows:
\begin{align*}
c^J_k &:= \#\{\text{$k$-cells of } J\} \\
m^J_k &:= \#\{\text{critical $k$-cells for $f|_J$}\} \\
u^J_k &:= \sum_{\sigma^k \in J} \#\{\tau \in J \mid \sigma\overset{\mathrm{nc}}{\prec}\tau\} \\
d^J_k &:= \sum_{\sigma^k \in J} \#\{\nu \in J \mid \nu\overset{\mathrm{nc}}{\prec}\sigma\}
\end{align*}

If $f$ is Morse, then the following relations hold:
\begin{align*}
c^J_k &= m^J_k + u^J_k + d^J_k 
\hspace{20pt} \text{for any $k \in \{0,1,\dots,\dim \K\}$}
\\
u^J_k &= d^J_{k+1}
 \hspace{61pt} \text{for any $k \in \{0,1,\dots,\dim \K-1 \}$}
\end{align*}

Moreover, we set
\[
c_k := c^\K_k, \quad m_k := m^\K_k, \quad u_k := u^\K_k, \quad d_k := d^\K_k
\]
for any $k \in \{0,1,\dots,\dim \K\}$.
\end{definition}

\begin{proposition}\label{propmbpoly}
Let $f$ be a discrete Morse-Bott function on a finite CW complex $\K$. 
For any reduced collection $C\in\mathcal{C}_f$, there is a polynomial $r(t)$ with nonnegative integer coefficients such that
\[
\sum_{k=0}^{\dim \K} \# \{\sigma\in C\mid\dim\sigma=k\} t^k = P_t(C) + (1+t) r(t).
\]

Moreover, the following equality holds: 
\[
r(t) = \sum_{k=1}^{\dim \K} \bigl(\rank B_{k-1}^C \bigr) t^{k-1}
\]
\end{proposition}

\begin{proof}
For any integer $k \in \Z_{\ge0}$, the following two short exact sequences exist: 
\[
0 \rightarrow Z^C_k \rightarrow C_k(C) \rightarrow B_{k-1}^C \rightarrow 0
\]
\[
0 \rightarrow B_k^C \rightarrow Z^C_k \rightarrow Z_k^C / B_k^C \rightarrow 0
\]
Then we obtain the following equalities for any integer $k \in \Z_{\ge0}$: 
\[
\rank C_k(C) = \rank Z^C_k + \rank B_{k-1}^C
\]
\[
\rank Z^C_k = \rank B_k^C + \rank(Z_k^C / B_k^C)
\]
These imply the following equality: 
\begin{align*}
    \sum_{k=0}^{\dim\K} \# \{\sigma\in C\mid\dim\sigma=k\} t^k - P_t(C)
    &= \sum_{k=0}^{\dim\K} \{\rank C_k(C) - b^C_k\} t^k \\
    &= \sum_{k=0}^{\dim\K} (\rank B_{k-1}^C + \rank B_k^C) t^k \\
    &= (1+t)\sum_{k=1}^{\dim\K} (\rank B_{k-1}^C) t^{k-1}.
\end{align*}
By $\rank B_{k-1}^C \ge 0$ for any integer $k \in \Z_{>0}$, any coefficient of $r(t)$ is nonnegative.
\end{proof}

We state the following folklore result.

\begin{theorem}\label{thmmorsep}
Let $f$ be a discrete Morse function on a finite CW complex $\K$.
Then there is a polynomial $r(t)$ with nonnegative integer coefficients such that  
\begin{align}\label{eq:dmf_poincare_poly}
\sum_{k=0}^{\dim \K} m_k t^k = P_t(\K) + (1+t) r(t)
\end{align}
where $P_t(\K) = \sum_k b_k t^k$ is the Poincar{\'e} polynomial of $\K$ and
\[
r(t) = \sum_{k=1}^{\dim \K} \bigl(\rank B_{k-1} - d_k \bigr)t^{k-1}.
\]
\end{theorem}

Although this proof is similar to that of \cite[Theorem~3.36]{2004MorseIneq} in the continuous setting, to the best of the authors' knowledge, no direct proof has been given in the literature, so we include it here for completeness.

\begin{proof}
For any $k \in \Z_{\ge0}$, from the two short exact sequences
\[
0\rightarrow  Z_k  \rightarrow C_k(\K) \rightarrow  B_{k-1}  \rightarrow 0
\]
and
\[
0 \rightarrow  B_k \rightarrow Z_k \rightarrow  Z_k / B_k  \rightarrow 0, 
\]
we obtain the following equalities: 
\[
\rank C_k(\K) = \rank Z_k +\rank B_{k-1}  
\]
\[
\rank Z_k =\rank B_k +\rank( Z_k / B_k )=\rank B_k +  b_k
\]
We have the following equality: 
\begin{align*}
    \sum_{k=0}^{\dim\K}m_k t^k-P_t(\K)&=\sum_{k=0}^{\dim\K}\{\rank C_k(\K)-u_k-d_k-b_k\}t^k\\&=\sum_{k=0}^{\dim\K}(\rank B_{k-1} -d_k+\rank B_k -d_{k+1})t^k\\&=(1+t)\sum_{k=1}^{\dim\K}(\rank B_{k-1} -d_k)t^{k-1}.
\end{align*}
In the sequence of equalities above, the relations $c_k = m_k + u_k + d_k$, $c_k = \operatorname{rank} C_k(\K)$ for any $k \in \Z_{\geq 0}$, and the definition of $P_t(\K)$ are applied in the first line. 
In the second line, the equality $b_k = \operatorname{rank}(Z_k/B_k)$ together with the two equalities above and the relation $d_{k+1} = u_k$ are used for any $k \in \Z_{\geq 0}$.

For a subcomplex \(J\) of \(\K\), for any $k \in \Z_{\geq 0}$, denote by $\partial_{J,k}: C_k(J) \to C_{k-1}(J)$ the boundary homomorphism.
For any $k \in \Z_{> 0}$, since $\rank B_{k-1} - d_k = u_k + m_k - \rank Z_k$, it suffices to show the inequality
\begin{align}\label{eq:non-neg}
u_k^J + m_k^J - \rank \ker(\partial_{J,k}) \;\geq\; 0 
\end{align}
for any subcomplex \(J\).
Assume that there is a subcomplex not satisfying the inequality~$(\ref{eq:non-neg})$.
Fix a minimal subcomplex \(K\) not satisfying the inequality~$(\ref{eq:non-neg})$.
Here, a minimal subcomplex not satisfying the inequality~$(\ref{eq:non-neg})$ means a subcomplex whose number of cells is smallest in the set of such subcomplexes not satisfying the inequality~$(\ref{eq:non-neg})$.
Note that $J = \{\sigma^0\}$ satisfies the inequality~$(\ref{eq:non-neg})$.
Let \(\sigma \in K\) be one of the maximal-dimensional cells of \(K\) such that
\[
f(\sigma) = \max \{ f(\alpha) \mid \alpha \text{ is a maximal-dimensional cell in } K \}.
\]
Then \(K_\sigma := K \setminus \{\sigma\}\) is a subcomplex of \(\K\) . 
Put $p := \dim \sigma$. 
Then $u^K_p = u^{K_\sigma}_p = 0$.
Since $K$ is a minimal subcomplex which does not satisfy $(\ref{eq:non-neg})$, the subcomplex $K_\sigma$ satisfies $(\ref{eq:non-neg})$. 
For any $k \in \{0,1, \ldots , \dim \K \}$, we have the following inequality: 
\[
\begin{split}
u_k^{K_\sigma} + m_k^{K_\sigma} - \rank \ker(\partial_{K_\sigma,k}) \geq 0 
\end{split}
\]
Because $K_\sigma$ is obtained from $K$ by removing the single cell $\sigma$, we obtain the following relations: 
\[
\begin{split}
u_p^{K} + m_p^{K} - (u_p^{K_\sigma} + m_p^{K_\sigma})  \in \{0, 1\} 
\\
\rank\ker\partial_{K,p} - \rank\ker\partial_{K_\sigma,p} \in \{0,1\}
\end{split}
\]
\begin{align}
\rank\ker\partial_{K,k} = \rank\ker\partial_{K_\sigma,k}
\label{eq:001a}
\end{align}
for any $k \neq p \in \{0,1, \ldots , \dim \K \}$. 
By the inequality $(\ref{eq:non-neg})$, we have the following inequalities: 
\[
u_p^K + m_p^K - \rank \ker(\partial_{K,p}) = -1 
\]
\[
u_p^{K_\sigma} + m_p^{K_\sigma} - \rank \ker(\partial_{K_\sigma,p}) = 0 
\]

Therefore, the following equalities hold:
\begin{align}
&u^K_p + m^K_p = u^{K_\sigma}_p + m^{K_\sigma}_p, \label{eq:001}\\
&\rank \ker \partial_{K,p} = \rank \ker \partial_{K_\sigma,p} + 1. \label{eq:002}
\end{align}

From $(\ref{eq:001a})$ and $(\ref{eq:002})$, there is no maximal-dimensional cell whose dimension is not \(p\).
Moreover, by $(\ref{eq:002})$, there is a $p$-cycle 
\[
x := \sum_{\alpha^p \in K} a_\alpha \alpha \in \ker \partial_{K,p}, \qquad a_\alpha \in \Z,
\]
with \(a_\sigma \ne 0\).
By $(\ref{eq:001})$, the cell \(\sigma\) is downward noncritical, so there is \(\nu \in K\) with \(\nu \overset{\mathrm{nc}}{\prec} \sigma\).
Since \(\partial_{K,k} x = 0\), there is \(\widetilde{\sigma} \in K \setminus \{ \sigma \} = K_\sigma \) with $a_{\widetilde{\sigma}} \ne0$ and \(\widetilde{\sigma} \succ \nu\) such that $f(\widetilde{\sigma}) > f(\nu)$. 
Then $f(\widetilde{\sigma}) > f(\nu) \ge f(\sigma)$.
Since \(\dim \widetilde{\sigma} \geq p\), the cell \(\widetilde{\sigma}\) is $p$-dimensional in \(K\), which contradicts the choice of \(\sigma\). 
\end{proof}

For any $k \in \Z_{\geq 0}$, denote by $\K^k$ the set of $k$-cells of $\K$.

\begin{lemma}\label{lemmbpoly}
Let $f$ be a discrete Morse-Bott function on a finite CW complex $\K$,  $k \in \{0,1,\dots,\dim \K\}$ an integer,
and  $x \in Z_k \setminus \{0\}$.
Write
\[
x = \sum_{\sigma \in \K^k} a_\sigma \,\sigma .
\]
and define a non-empty set $A := \{\sigma \in \K^k \mid a_\sigma \neq 0\}$.
For any cell \(\alpha \in A\) with
\begin{align}\label{eq:maxalpha}
    f(\alpha) = \max\{\, f(\sigma) \mid \sigma \in A \,\}
\end{align}
the sum
\[
x^{\alpha} := \sum_{\sigma \in A \cap \mathcal{C}_f(\alpha)} a_\sigma \sigma
\]
satisfies \(x^\alpha \in Z^{\mathcal{C}_f(\alpha)}_k \setminus \{0\}\).
\end{lemma}

\begin{proof}
Put $C := \mathcal{C}_f(\alpha)$. 
First, we show the following statement.  

\begin{claim}\label{claimempty}
\(\{\nu \prec \sigma \mid \nu \in C\} = \emptyset\) for any cell \(\sigma \in A \setminus C\).
\end{claim}

\begin{proof}
Fix any cell \(\sigma \in A \setminus C\). 
Define $N_\sigma := \{\nu \prec \sigma \mid \nu \in C\}$. 
Assume that \(N_\sigma \neq \emptyset\).  
Fix any cell \(\nu \in N_\sigma\).  
Since \(\nu \in C\), the cell $\nu$ is not upward noncritical and so
\[
f(\alpha) = f(\nu) \le f(\sigma).
\]
Because \(\sigma \in A\), by the maximality (\ref{eq:maxalpha}) of $\alpha$, it follows \(f(\nu) = f(\sigma)\).  
If \(\sigma\) is upward noncritical, then by Lemma~\ref{lemupnon}, \(\nu\) is also upward noncritical, which contradicts that \(\nu\) is not.  
Thus \(\sigma\) is downward noncritical.
Then there is \(\widetilde{\nu} \prec \sigma\) such that \(f(\widetilde{\nu}) > f(\sigma)\).  
By \(\partial_k x = 0\), there is \(\widetilde{\sigma} \in A \setminus \{ \sigma \} \) with \(\widetilde{\nu} \prec \widetilde{\sigma}\).
Then $f(\alpha) = f(\sigma) < f(\widetilde{\nu}) \le f(\widetilde{\sigma})$, which contradicts the maximality (\ref{eq:maxalpha}) of $\alpha$. 
\end{proof}

By the previous claim, we have the following relation: 
\begin{align*}
\{\nu \in \K \mid \sigma \in A, \; \nu \prec \sigma\}
&= \{\nu \in \K \mid \sigma \in A, \; \nu \prec \sigma, \; \nu \notin C\} \\
&\quad \cup \{\nu \in \K \mid \sigma \in A \cap C, \; \nu \prec \sigma, \; \nu \in C\}.
\end{align*}

Therefore, we have the following equality: 
\begin{align*}
0 = \partial_k x   &= \sum_{\sigma \in A} \sum_{\nu \prec \sigma} a_\sigma [\sigma:\nu] \nu \\
  &= \sum_{\sigma \in A \cap C} \sum_{\substack{\nu \prec \sigma \\ \nu \in C}} a_\sigma [\sigma:\nu] \nu
     + \sum_{\sigma \in A} \sum_{\substack{\nu \prec \sigma \\ \nu \notin C}} a_\sigma [\sigma:\nu] \nu
\end{align*}

By the independence of the basis of the chain complex, each term in the above equation vanishes. 
In particular, we have the following equality: 
    
\begin{align*}     
  0 &= \sum_{\sigma \in A \cap C} \sum_{\substack{\nu \prec \sigma \\ \nu \in C}} a_\sigma [\sigma:\nu] \nu
   = \partial^C_k x^\alpha
\end{align*}
\end{proof}

We also introduce the following notation to clarify the following discussion.

\begin{definition}
Let $f:\K\to\R$ be a function, $k \in \{0,1,\dots,\dim \K\}$ an integer, and $J$ a subcomplex of $\K$. 
Then $\tilde{m}^J_k,\tilde{u}^J_k,\tilde{d}^J_k$ are defined as follows: 
\begin{align*}
\tilde{m}^J_k &:= \#\{\text{weakly critical $k$-cells for $f|_J$}\} \\
\tilde{u}^J_k &:= \sum_{\sigma^k \in J} \#\{\tau \in J \mid \sigma\overset{\mathrm{snc}}{\prec}\tau\} \\
\tilde{d}^J_k &:= \sum_{\sigma^k \in J} \#\{\nu \in J \mid \nu\overset{\mathrm{snc}}{\prec}\sigma\}
\end{align*}

If $f$ is Morse-Bott, then the following relations hold:
\begin{align*}
c^J_k &= \tilde{m}^J_k + \tilde{u}^J_k + \tilde{d}^J_k 
\hspace{20pt} 
\text{ for any $k \in \{0,1,\dots,\dim \K\}$}
\\
\tilde{u}^J_k &= \tilde{d}^J_{k+1}
\hspace{61pt} 
\text{ for any $k \in \{0,1,\dots,\dim \K -1 \}$}
\end{align*}

Moreover, we set
\[
\tilde{m}_k := \tilde{m}^\K_k, \quad \tilde{u}_k := \tilde{u}^\K_k, \quad \tilde{d}_k := \tilde{d}^\K_k.
\]
for any $k \in \{0,1,\dots,\dim \K\}$.
\end{definition}

Now, we show the following inequality. 

\begin{theorem}\label{thmmbpoly}
Let $f$ be a discrete Morse-Bott function on a finite CW complex $\K$. 
Then there is a polynomial $R(t)$ with nonnegative integer coefficients such that
\begin{align}\label{eq:dmbf_poincare_poly}
\sum_{C \in \mathcal{C}_f} P_t(C) \;=\; P_t(\K) + (1+t)R(t)    
\end{align}
where $R(t) := \sum_{k=1}^{\dim \K} \Bigl(\rank B_{k-1} - \tilde{d}_k - \sum_{C \in \mathcal{C}_f}\rank B_{k-1}^C\Bigr)t^{k-1}$.
\end{theorem}

\begin{proof}
Since $\K$ is finite, the combinatorial vector field $-\nabla_s f$ is positively bounded.  
By Lemma~\ref{lemsgvf}, there is a discrete Morse function $g$ such that $-\nabla_s f = -\nabla g$.  
Applying Theorem~\ref{thmmorsep} to $g$, we obtain the following equality:
\[
\sum_{k=0}^{\dim \K} m^g_k t^k = P_t(\K) + (1+t) r(t),
\]
where 
\[
r(t) = \sum_{k=1}^{\dim \K} \bigl(\rank B_{k-1} - d^g_k\bigr) t^{k-1}.
\]
The relation $-\nabla_s f = -\nabla g$ implies $d^g_k = \tilde{d}_k$ and $u^g_k = \tilde{u}_k$. 
Then $m^g_k = \tilde{m}_k$.
Therefore, the following holds:
\begin{align*}
    &\sum_{k=0}^{\dim\K}m^g_k t^k=P_t(\K)+(1+t)r(t)\\
\implies & \sum_{k=0}^{\dim\K}\tilde{m}_k t^k=P_t(\K)+(1+t)r(t)\\
\implies & \sum_{k=0}^{\dim\K}\sum_{C\in\mathcal{C}_f}\# \{\sigma\in C\mid\dim\sigma=k\} t^k=P_t(\K)+(1+t)r(t)\\
\implies & \sum_{C\in\mathcal{C}_f}(P_t(C)+(1+t)\sum_{k=1}^{\dim\K}(\rank B_{k-1}^C)t^{k-1})=P_t(\K)+(1+t)r(t)\\
\implies & \sum_{C\in\mathcal{C}_f}P_t(C)=P_t(\K)+(1+t)\\
& \hspace{110pt} \sum_{k=1}^{\dim\K}(\rank B_{k-1} -d^g_k-\sum_{C\in\mathcal{C}_f}\rank B_{k-1}^C)t^{k-1}.
\end{align*}

In the sequence of equalities above, the definition of $\tilde{m}_k$ is applied in the third line, Proposition~\ref{propmbpoly} in the fourth, and the definition of $r(t)$ in the fifth.

It suffices to show that the coefficients of $R(t)$ are nonnegative.
Fix any $k \in \{1,2,\dots,\dim \K\}$. 
The following equalities 
\begin{align*}
\rank C_k(\K) &= \rank Z_k +\rank B_{k-1} \\
\rank C_k(C) &= \rank Z^C_k + \rank B_{k-1}^C\\
c^J_k &=\tilde{m}^J_k + \tilde{u}^J_k + \tilde{d}^J_k
\end{align*}
imply the following equality:
\[
\rank B_{k-1} - \tilde{d}_k - \sum_{C \in \mathcal{C}_f} \rank B_{k-1}^C
= \tilde{u}_k + \sum_{C \in \mathcal{C}_f} \rank Z^C_k - \rank Z_k.
\]
For any subcomplex $J$ of $\K$, put $\mathcal{C}_J := \mathcal{C}_{f\mid_J}$. 
Then it suffices to show the following inequality:
\begin{align}\label{eq:ineq}
\tilde{u}^J_k + \sum_{C \in \mathcal{C}_J} \rank Z^C_k - \rank \ker \partial_{J,k} \geq 0.
\end{align}

Assume that there is a subcomplex of $\K$ that does not satisfy $(\ref{eq:ineq})$.  
Note that $J := \{\sigma^0\}$ satisfies $(\ref{eq:ineq})$.  
Hence there is a minimal subcomplex of $\K$ not satisfying $(\ref{eq:ineq})$, and denote it by $K$.  
Define
\[
M := \{\sigma \in K \mid \sigma \text{ is maximal-dimensional in } K\}.
\]

For any $\sigma \in M$, put $K_\sigma := K \setminus \{\sigma\}$.  
Then $K_\sigma$ is a subcomplex of $\K$ and satisfies $(\ref{eq:ineq})$.  
Therefore, by the same argument in the proof of Theorem~\ref{thmmorsep}, the following two relations hold for any $\sigma \in M$:
\begin{align}
&\tilde{u}^K_k + \sum_{C \in \mathcal{C}_K} \rank Z^C_k
= \tilde{u}^{K_\sigma}_k + \sum_{C \in \mathcal{C}_{K_\sigma}} \rank Z^C_k, \label{eq:03}\\
&\rank \ker \partial_{K,k}
= \rank \ker \partial_{K_\sigma,k} + 1. \label{eq:04}
\end{align}

From $(\ref{eq:04})$ it follows that $\dim \sigma = k$ for any $\sigma \in M$.  
Therefore, all elements of $M$ are $k$-cells, and conversely, every $k$-cell belongs to $M$. 
Choose a $k$-cell $\alpha \in M$ such that $f(\alpha) = \max \{ f(\sigma) \mid \sigma \in M\}$.
By $(\ref{eq:04})$, there is a chain $x \in \ker \partial_{K,k}$ which can be written as
\[
x = \sum_{\sigma^k \in K} a_\sigma \sigma
\]
with $a_\alpha \ne 0$.  
Assume that $\alpha$ is downward noncritical.  
Then there is a cell $\nu \in K$ with $\nu \overset{\mathrm{snc}}{\prec} \alpha$.  
Since $\partial_{K,k} x = 0$, there is a maximal-dimensional cell $\widetilde{\sigma} \succ \nu$ in $K$ with $\widetilde{\sigma} \neq \alpha$ such that $f(\alpha) < f(\nu) \le f(\widetilde{\sigma})$, which contradicts the maximality of $\alpha$.

Thus, put $C := \mathcal{C}_K(\alpha)$.
If the set difference $C\setminus\{\alpha\}$ is a reduced collection, then Lemma~\ref{lemmbpoly} implies $\rank Z_k^C\ne\rank Z_k^{C\setminus\{\alpha}\}$, which contradicts $(\ref{eq:03})$.
Thus, the set difference $C\setminus\{\alpha\}$ is a finite disjoint union of reduced collections $C_1 \sqcup \cdots \sqcup C_l$. 
By Lemma~\ref{lemmbpoly}, we have $\rank Z_k^C\ne \sum_{i=1}^l \rank Z_k^{C_i}$, which contradicts $(\ref{eq:03})$.
\end{proof}

We show that the equality (\ref{eq:dmbf_poincare_poly})  is a generalization of the equality (\ref{eq:dmf_poincare_poly}) as follows. 

\begin{proposition}\label{rmkmb}
Let $f$ be a discrete Morse function on a finite CW complex $\K$.
Then the equality (\ref{eq:dmbf_poincare_poly}) is reduced into the equality (\ref{eq:dmf_poincare_poly}).
\end{proposition}

\begin{proof}
By definition of regularity of faces, we obtain that  $P_t(C) = 0$ for any reduced collection $C \in \mathcal{C}_f$ with  $C = \{\sigma, \tau\}$ with $\sigma^p \overset{\mathrm{reg}}{\prec} \tau$. 
Suppose that $f$ is a discrete Morse function. By Corollary~\ref{cormorse-mb}, for any reduced collection $C\in\mathcal{C}_f$, we obtain that $\#C\le2$ and $C=L(C)$.

\begin{claim}
$\sum_{C \in \mathcal{C}_f} P_t(C)
= \sum_{k=0}^{\dim \K} m_k t^k$.
\end{claim}

\begin{proof}
For any reduced collection $C\in\mathcal{C}_f$ with $\# C=2$, we have that $P_t(C)=0$. 
For any reduced collection $C\in\mathcal{C}_f$ with $\# C=1$, by Remark~\ref{rmk:3}, $\#C=\#L(C)=1$ implies $C=L(C)=\{\text{a critical cell of } f\vert_C \}$ and so $P_t(C) = t^{\dim C}$.
Therefore, the desired equality holds. 
\end{proof}

\begin{claim}
$d_k=\tilde{d_k}+\sum_{C \in \mathcal{C}_f}\rank B_{k-1}^C$.
\end{claim}

\begin{proof}
For any reduced collection $C \in \mathcal{C}_f$ with $\# C = 1$, the definition of $B_{k-1}^C$ implies $\rank B_{k-1}^C =0$. For any reduced collection $C \in \mathcal{C}_f$ with $\#C=2$, we have $C=\{\nu\prec\sigma^p\}$ and so $\rank B_{p-1}^C=\rank \im\partial_p^C=1$. By definition of $d_k$, we have the following equalities:
\begin{align*}
    d_k&=\sum_{\sigma^k\in\K}\#\{\nu\in\K\mid\nu\overset{\mathrm{nc}}{\prec}\sigma\}\\
    &=\sum_{\sigma^k\in\K}\#\{\nu\prec\sigma\mid f(\nu)>f(\sigma)\}+\sum_{\sigma^k\in\K}\#\{\nu\prec\sigma\mid f(\nu)=f(\sigma)\}\\
    &=\tilde{d_k}+\#\{C\in\mathcal{C}_f\mid C=\{\nu\prec\sigma^k\},~\nu,\sigma\in\K\}\\
    &=\tilde{d_k}+\sum_{C\in\mathcal{C}_f}\rank B_{k-1}^C.
\end{align*}
In the sequence of equalities above, by $\#C\le2$, the third equality holds.
\end{proof}

The previous claims imply the assertion. 
\end{proof}

\section*{Acknowledgements}
The authors would like to thank Professor Jesus Gonzalez for his helpful comments. 

\bibliographystyle{abbrv}
\bibliography{DMBF}

\end{document}